\documentclass[11pt,letterpaper,reqno]{amsart}
\usepackage{epsfig}

\def\be{\begin{equation}}
\def\ee{\end{equation}}
\def\bea{\begin{eqnarray}}
\def\eea{\end{eqnarray}}
\def\bes{\begin{eqnarray*}}
\def\ees{\end{eqnarray*}}

\def\nn{\nonumber}
\def\lb{\label}
\def\bs{\setminus}

\def\R{{\bf R}}
\def\C{{\bf C}}
\def\Z{{\bf Z}}

\def\N{{\bf N}}
\def\U{{\bf U}}

\def\Q{{\bf Q}}

\def\aa{{\alpha}}
\def\bb{{\beta}}
\def\ga{{\gamma}}

\def\th{{\theta}}

\def\om{{\omega}}
\def\Om{{\Omega}}

\def\lm{{\lambda}}
\def\Lm{{\Lambda}}
\def\dl{{\delta}}

\def\sg{{\sigma}}

\def\<{{\langle}}
\def\>{{\rangle}}

\def\rank{{\rm rank}}

\def\Sp{{\rm Sp}}

\def\ol{\overline}

\def\hb{\vrule height0.18cm width0.14cm $\,$}

\title{Multiple closed geodesics on Finsler $3$-dimensional sphere}

\author[Huagui Duan]{Huagui Duan}
\thanks{Huagui Duan was partially supported by National Key R\&D Program of China (Grant No. 2020YFA0713300)
and NNSFC (Nos. 12271268, 11671215 and 12361141812) and Nankai University. Zihao Qi was partially supported by NNSFC (No. 12271268).}
\address{Huagui Duan, School of Mathematical Sciences and LPMC, Nankai University, Tianjin 300071, The People's Republic of China}
\email{duanhg@nankai.edu.cn.}

\author[Zihao Qi]{Zihao Qi}
\address{Zihao Qi, School of Mathematical Sciences, Nankai University, Tianjin 300071, The People's Republic of China}
\email{2120210039@mail.nankai.edu.cn.}

\begin{document}
\maketitle

\begin{abstract}
{\it In 1973, Katok constructed a non-degenerate (also called bumpy) Finsler metric on $S^3$ with exactly four prime closed geodesics. And then Anosov conjectured that four should be the optimal lower bound of the number of prime closed geodesics on every Finsler $S^3$. In this paper, we proved this conjecture for bumpy Finsler $S^{3}$ if the Morse index of any prime closed geodesic is nonzero.}
\end{abstract}



\renewcommand{\theequation}{\thesection.\arabic{equation}}
\renewcommand{\thefigure}{\thesection.\arabic{figure}}

\setcounter{figure}{0}
\setcounter{equation}{0}
\section{Introduction and main result}

A closed curve on a Finsler manifold is a closed geodesic if it is
locally the shortest path connecting any two nearby points on this
curve. As usual, on any Finsler manifold
$(M, F)$, a closed geodesic $c:S^1=\R/\Z\to M$ is {\it prime}
if it is not a multiple covering (i.e., iteration) of any other
closed geodesics. Here the $m$-th iteration $c^m$ of $c$ is defined
by $c^m(t)=c(mt)$. The inverse curve $c^{-1}$ of $c$ is defined by
$c^{-1}(t)=c(1-t)$ for $t\in \R$.  Note that unlike a Riemannian manifold,
the inverse curve $c^{-1}$ of a closed geodesic $c$
on a irreversible Finsler manifold need not be a geodesic.
We call two prime closed geodesics
$c$ and $d$ {\it distinct} if there is no $\th\in (0,1)$ such that
$c(t)=d(t+\th)$ for all $t\in\R$.
On a reversible Finsler (or Riemannian) manifold, two closed geodesics
$c$ and $d$ are called { \it geometrically distinct} if $
c(S^1)\neq d(S^1)$, i.e., their image sets in $M$ are distinct.
We shall omit the word {\it distinct} when we talk about more than one prime closed geodesic.

For a closed geodesic $c$ on $n$-dimensional manifold $(M,\,F)$, denote by $P_c$
the linearized Poincar\'{e} map of $c$. Then $P_c\in \Sp(2n-2)$ is symplectic.
For any $P\in \Sp(2k)$, we define the {\it elliptic height } $e(P)$
of $P$ to be the total algebraic multiplicity of all eigenvalues of
$P$ on the unit circle $\U=\{z\in\C|\; |z|=1\}$ in the complex plane
$\C$. Since $P$ is symplectic, $e(P)$ is even and $0\le e(P)\le 2k$.
A closed geodesic $c$ is called {\it elliptic} if $e(P_c)=2(n-1)$, i.e., all the
eigenvalues of $P_c$ locate on $\U$; {\it non-degenerate} if $1$ is not an eigenvalue of $P_c$. A Finsler manifold $(M,\,F)$
is called {\it bumpy} if all the closed geodesics (including iterations) on it are non-degenerate.

There is a famous conjecture in Riemannian geometry which claims there exist infinitely many
closed geodesics on any compact Riemannian manifold. This conjecture
has been proved except for CROSS's (compact rank one symmetric
spaces). The results of Franks \cite{Fra} in 1992, Hingston \cite{Hin2} and Bangert \cite{Ban}
in 1993 imply this conjecture is true for any Riemannian $S^2$.

For the Finsler metric, the closed geodesic problem is completely different. It was quite surprising when Katok \cite{Kat}  in 1973 found
some irreversible bumpy Finsler metrics on CROSS's with only finitely many closed geodesics and all closed geodesics are non-degenerate and elliptic (cf. \cite{Zil}). Based on Katok's examples, in 1974, Anosov \cite{Ano} conjectured that the optimal lower bound of the number of prime closed geodesics on any Finsler $S^n$ should be $2[\frac{n+1}{2}]$. We refer readers to two excellent survey papers \cite{Lon4} and \cite{BuK} for more interesting questions.

The index iteration theory of closed geodesics (cf.\cite{Bot} and \cite{Lon3}) has been an important and powerful tool in studying the closed geodesic problem on Finsler manifolds. For example, Bangert and Long in \cite{BaL} (finished in 2005) shows that there exist at least two prime closed geodesics on every Finsler $S^2$, which solved the Anosov conjecture. And then a great number of multiplicity and stability results about closed geodesics which are established by index iteration theory have appeared (cf. \cite{DuL1}-\cite{DuL2}, \cite{DLW}, \cite{Lon4}, \cite{Rad3}, \cite{Wan1}-\cite{Wan3} and therein).

This paper is mainly devoted to study the Anosov conjecture of bumpy Finsler $S^3$. In this direction, Duan and Long in \cite{DuL1} and Rademacher in \cite{Rad3} proved the existence of at least two prime closed geodesics on bumpy Finsler $S^n$, respectively. Wang in \cite{Wan2} solved the Anosov conjecture for any positively curved bumpy Finsler $S^n$. Recently, the Anosov conjecture on bumpy Finsler CROSS's (including spheres) has been  proved by Duan, Long and Wang in \cite{DLW} under the much weaker curvature or index conditions, and especially for $S^3$ the index condition is $i(c)\ge 2$ for any prime closed geodesic $c$, which, to the authors' knowledge, is the weakest condition for Anosov conjecture of bumpy Finsler $S^3$. Actually we also notice that in Katok's emxample, the Morse index of every prime closed geodesic on $S^3$ is equal to and greater than $2$.

Motivated by the Katok's examples and above results, in this paper we established the following result, which solved the Anosov conjecture of bumpy $S^3$ under much more relaxing restriction.

\vspace{2mm}

{\bf Theorem 1.1.} {\it Assume that $F$ is a bumpy Finsler metric on $S^{3}$, then there exist at least four prime closed geodesics on $(S^3,F)$ if the Morse index of every prime closed geodesic is nonzero.}

\vspace{2mm}

Under the assumption of Theorem 1.1, it follows from Remark 3.7 of \cite{DuL2} and Theorem 1.4 of \cite{DLW} that there exist at least three prime closed geodesics on bumpy Finsler $S^3$. Then, through assuming the existence of exactly three prime closed geodesics, we will use a contradiction argument to complete the proof of Theorem 1.1.

Here notice that we only assume the Morse index $i(c)\ge 1$ for any prime closed geodesic $c$, in general the sequence $\{i(c^m)\}$ will have no monotonicity with respect to $m$, which will lead to the invalidness of methods in \cite{DLW} and \cite{Wan2}. On the other hand, due to the assumption of the existence of three prime closed geodesics, the methods in \cite{DuL2} through some complicated classification arguments and the mean index identity can not give any valuable information in this situation.

To overcome the above difficulties, the main novelty in this paper is to establish a {\it variant of generalized common index jump theorem} (see Theorem 3.8 below), which, compared with Theorem 3.6 of \cite{DLLW} (cf. Theorem 3.6 below), gives more information about some parameters including the common integer $N$ and iterations $m_k$s. Then we use the Morse theory, some precise analysis and estimates of common iteration indices to derive a contradiction. Because this method does not use the special classification of closed geodesics, it may be hopeful to be applied to studying the corresponding problems on higher Finsler manifolds.

When the assumption of Theorem 1.1 is not satisfied, i.e. such $(S^3,F)$ possesses at least one prime closed geodesic with zero Morse index, the existence of infinitely many prime closed geodesics is suspected by some recent discussions of the first author, Y. Long and C. Zeng. Note that when a closed geodesic $c$ is non-degenerate, then the non-local minimum of $c$ with respect to the energy functional $E$ implies that the Morse index of $c$ satisfies $i(c)>0$. In fact, they proposed the following conjecture.

\vspace{2mm}

{\bf Conjecure 1.2.} {\it Let $(M,F)$ be a compact simply connected $n$-dimensional irreversible Finsler manifold. If there exist only finitely many prime closed geodesics on $(M,F)$, then none of the closed geodesics including all the iterates is a local minimum of the energy functional $E$.}

\vspace{2mm}

A related result with this is the Corollary 2 of \cite{BaK}, which shows that for a compact manifold $M$ with finite fundamental group, suppose there exists a closed geodesic $c$ on $M$ such that $c^m$ is a local minimum of $E$ for infinitely many $m\in\N$, then there exist infinitely many prime closed geodesics on $M$.

In this paper, let $\N$, $\N_0$, $\Z$, $\Q$, $\R$, and $\C$ denote the sets of natural integers, non-negative integers, integers,
rational numbers, real numbers, and complex numbers respectively. We use only singular homology modules with $\Q$-coefficients.
For an $S^1$-space $X$, we denote by $\overline{X}$ the quotient space $X/S^1$. We define the functions
\be \left\{\begin{array}{ll}[a]=\max\{k\in\Z\,|\,k\le a\},\quad
           E(a)=\min\{k\in\Z\,|\,k\ge a\} , \\
    \varphi(a)=E(a)-[a],\quad \{a\}=a-[a]. \lb{1.1}
    \end{array}\right.\ee

Especially, $\varphi(a)=0$ if $ a\in\Z\,$, and $\varphi(a)=1$ if $
a\notin\Z\,$.

\setcounter{equation}{0}

\section{Critical modules of closed geodesics and Morse Theory}

Let $M=(M,F)$ be a compact Finsler manifold $(M,F)$, the space
$\Lambda=\Lambda M$ of $H^1$-maps $\gamma:S^1\rightarrow M$ has a
natural structure of Riemannian Hilbert manifold on which the
group $S^1=\R/\Z$ acts continuously by isometries, cf.
\cite{Kli}, Chapters 1 and 2. This action is defined by
$(s\cdot\gamma)(t)=\gamma(t+s)$ for all $\gamma\in\Lm$ and $s,
t\in S^1$. For any $\gamma\in\Lambda$, the energy functional is
defined by
\be E(\gamma)=\frac{1}{2}\int_{S^1}F(\gamma(t),\dot{\gamma}(t))^2dt.
\lb{2.1}\ee
It is $C^{1,1}$ and invariant under the $S^1$-action. The
critical points of $E$ of positive energies are precisely the closed geodesics
$\gamma:S^1\to M$. The index form of the functional $E$ is well
defined along any closed geodesic $c$ on $M$, which we denote by
$E''(c)$. As usual, we denote by $i(c)$ and
$\nu(c)$ the Morse index and nullity of $E$ at $c$. In the
following, we denote by
\be \Lm^\kappa=\{d\in \Lm\;|\;E(d)\le\kappa\},\quad \Lm^{\kappa-}=\{d\in \Lm\;|\; E(d)<\kappa\},
  \quad \forall \kappa\ge 0. \nn\ee
For a closed geodesic $c$ we set $ \Lm(c)=\{\ga\in\Lm\mid E(\ga)<E(c)\}$.

For $m\in\N$ we denote the $m$-fold iteration map
$\phi_m:\Lambda\rightarrow\Lambda$ by $\phi_m(\ga)(t)=\ga(mt)$, for all
$\,\ga\in\Lm, t\in S^1$, as well as $\ga^m=\phi_m(\gamma)$. If $\gamma\in\Lambda$
is not constant then the multiplicity $m(\gamma)$ of $\gamma$ is the order of the
isotropy group $\{s\in S^1\mid s\cdot\gamma=\gamma\}$. For a closed geodesic $c$,
the mean index $\hat{i}(c)$ is defined as usual by
$\hat{i}(c)=\lim_{m\to\infty}i(c^m)/m$. Using singular homology with rational
coefficients we consider the following critical $\Q$-module of a closed geodesic
$c\in\Lambda$:
\be \overline{C}_*(E,c)
   = H_*\left((\Lm(c)\cup S^1\cdot c)/S^1,\Lm(c)/S^1\right). \lb{2.3}\ee

The following results of Rademacher will be used in our proofs below.

{\bf Proposition 2.1.} (cf. Satz 6.11 of \cite{Rad2} ) {\it Let $c$ be a
prime closed geodesic on a bumpy Finsler manifold $(M,F)$. Then there holds}
$$ \overline{C}_q( E,c^m) = \left\{\begin{array}{ll}
     \Q, \quad {\it if}\;\; i(c^m)-i(c)\in 2\Z\;\;{\it and}\;\;
                   q=i(c^m),\\
     0, \quad {\it otherwise}. \end{array}\right.  $$

{\bf Definition 2.2.} (cf. Definition 1.6 of \cite{Rad1}) {\it For a
closed geodesic $c$, let $\ga_c\in\{\pm\frac{1}{2},\pm1\}$ be the
invariant defined by $\ga_c>0$ if and only if $i(c)$ is even, and
$|\ga_c|=1$ if and only if $i(c^2)-i(c)$ is even. }

Let $(X,Y)$ be a space pair such that the Betti numbers
$b_i=b_i(X,Y)=\dim H_i(X,Y;\Q)$ are finite for all $i\in \Z$. As
usual the Poincar\'e series of $(X,Y)$ is defined by the
formal power series $P(X, Y)=\sum_{i=0}^{\infty}b_it^i$. We need
the following well known results on Betti numbers and the Morse
inequality for $\overline{\Lm}\equiv
\overline{\Lm} S^3$ and $\ol{\Lm}^0=\ol{\Lambda}^0S^3
=\{{\rm constant\;point\;curves\;in\;}S^3\}\cong S^3$.

{\bf Proposition 2.3.} (cf. Remark 2.5 of \cite{Rad1} or \cite{Hin1}) {\it The
Poincar\'e series is given by \bea P(\ol{\Lm}S^3,\ol{\Lm}^0S^3)(t)
&=&t^2\left(\frac{1}{1-t^2}+\frac{t^2}{1- t^2}\right)  \nn\\
&=& t^2(1+t^2)(1+t^2+t^4+\cdots) = t^2+2t^4+2t^6+\cdots, \nn\eea
which yields}
\bea {b}_q &=& {b}_q(\ol{\Lm}S^3,\ol{\Lm}^0 S^3)\;
 = \rank H_q(\ol{\Lm} S^3,\ol{\Lm}^0 S^3 )\nn\\
 &=& \;\left\{\begin{array}{lll}
    1,\quad {\it if}\quad q=2,  \\
    2,\quad {\it if}\quad q=2k+2,\quad k\in \N, \\
    0 \quad {\it otherwise}. \end{array}\right. \lb{2.4}\eea

{\bf Proposition 2.4.} (cf. Theorem I.4.3 of \cite{Cha}, Theorem
6.1 of \cite{Rad2}) {\it Suppose that there exist only finitely
many prime closed geodesics $\{c_j\}_{1\le j\le k}$ on a Finsler
$3$-sphere $(S^3, F)$. Set
$$ M_q =\sum_{1\le j\le k,\; m\ge 1}\dim{\ol{C}}_q(E, c^m_j), \quad \forall q\in\Z. $$
Then for every integer $q\ge 0$ there holds }
\bea
M_q - M_{q-1} + \cdots +(-1)^{q}M_0
&\ge& b_q - b_{q-1}+ \cdots + (-1)^{q}b_0, \lb{2.5}\\
M_q &\ge& b_q. \lb{2.6}\eea

\setcounter{figure}{0}
\setcounter{equation}{0}
\section{The common index jump theorem for symplectic paths}

In \cite{Lon1} of 1999, Y. Long established the basic normal form decomposition of symplectic matrices. Based on this result he
further established the precise iteration formulae of indices of symplectic paths in \cite{Lon2} of 2000. Since every closed geodesic on a sphere is orientable, then, by Theorem 1.1 of \cite{Liu}, the Morse index of a closed geodesic on $S^n$ coincides with the Maslov-type index of a corresponding symplectic path.

As in \cite{Lon2}, the basic normal forms are denoted by
\bea
N_1(\lm, b) &=& \left(\begin{array}{ll}\lm & b\\
                                0 & \lm \end{array}\right), \qquad {\rm for\;}\lm=\pm 1, \; b\in\R, \lb{3.1}\\
D(\lm) &=& \left(\begin{array}{ll}\lm & 0\\
                      0 & \lm^{-1} \end{array}\right), \qquad {\rm for\;}\lm\in\R\bs\{0, \pm 1\}, \lb{3.2}\\
R(\th) &=& \left(\begin{array}{ll}\cos\th & -\sin\th \\
                           \sin\th & \cos\th \end{array}\right), \qquad {\rm for\;}\th\in (0,\pi)\cup (\pi,2\pi), \lb{3.3}\\
N_2(e^{\sqrt{-1}\th}, B) &=& \left(\begin{array}{ll} R(\th) & B \\
                  0 & R(\th) \end{array} \right), \qquad {\rm for\;}\th\in (0,\pi)\cup (\pi,2\pi)\;\; {\rm and}\; \nn\\
        && \quad B=\left(\begin{array}{lll} b_1 & b_2\\
                                  b_3 & b_4 \end{array}\right)\; {\rm with}\; b_j\in\R, \;\;
                                         {\rm and}\;\; b_2\not= b_3. \lb{3.4}\eea

As in \cite{Lon2}, the $\diamond$-sum (direct sum) of any two real matrices is defined by
$$ \left(\begin{array}{lll}A_1 & B_1\\ C_1 & D_1 \end{array}\right)_{2i\times 2i}\diamond
      \left(\begin{array}{lll}A_2 & B_2\\ C_2 & D_2 \end{array}\right)_{2j\times 2j}
=\left(\begin{array}{llll}A_1 & 0 & B_1 & 0 \\
                                   0 & A_2 & 0& B_2\\
                                   C_1 & 0 & D_1 & 0 \\
                                   0 & C_2 & 0 & D_2\end{array}\right). $$

For every $M\in\Sp(2n)$, the homotopy set $\Omega(M)$ of $M$ in $\Sp(2n)$ is defined by
$$ \Om(M)=\{N\in\Sp(2n)\,|\,\sg(N)\cap\U=\sg(M)\cap\U\equiv\Gamma,
                    \;\nu_{\om}(N)=\nu_{\om}(M),\, \forall\om\in\Gamma\}, $$
where $\sg(M)$ denotes the spectrum of $M$,
$\nu_{\om}(M)\equiv\dim_{\C}\ker_{\C}(M-\om I)$ for $\om\in\U$.
The component $\Om^0(M)$ of $P$ in $\Sp(2n)$ is defined by
the path connected component of $\Om(M)$ containing $M$.

\medskip

{\bf Lemma 3.1.} (cf. \cite{Lon2}, Lemma 9.1.5 and List 9.1.12 of \cite{Lon3})
{\it For $M\in\Sp(2n)$ and $\om\in\U$, the splitting number $S_M^\pm(\om)$ (cf. Definition 9.1.4 of \cite{Lon3}) satisfies
\begin{eqnarray}
S_M^{\pm}(\om) &=& 0, \qquad {\rm if}\;\;\om\not\in\sg(M).  \lb{3.5}\\
S_{N_1(1,a)}^+(1) &=& \left\{\begin{array}{lll}1, &\quad {\rm if}\;\; a\ge 0, \\
0, &\quad {\rm if}\;\; a< 0. \end{array}\right. \lb{3.6}\eea

For any $M_i\in\Sp(2n_i)$ with $i=0$ and $1$, there holds }
\be S^{\pm}_{M_0\diamond M_1}(\om) = S^{\pm}_{M_0}(\om) + S^{\pm}_{M_1}(\om),
    \qquad \forall\;\om\in\U. \lb{3.7}\ee

We have the following decomposition theorem

\medskip

{\bf Theorem 3.2.} (cf. \cite{Lon2} and Theorem 1.8.10 of \cite{Lon3}) {\it For
any $M\in\Sp(2n)$, there is a path $f:[0,1]\to\Om^0(M)$ such that $f(0)=M$ and
\be f(1) = M_1\diamond\cdots\diamond M_k,  \lb{3.8}\ee
where each $M_i$ is a basic normal form listed in (\ref{3.1})-(\ref{3.4})
for $1\leq i\leq k$.}

\medskip

For every $\ga\in\mathcal{P}_\tau(2n)\equiv\{\ga\in C([0,\tau],Sp(2n))\ |\ \ga(0)=I_{2n}\}$, we extend
$\ga(t)$ to $t\in [0,m\tau]$ for every $m\in\N$ by
\bea \ga^m(t)=\ga(t-j\tau)\ga(\tau)^j \quad \forall\;j\tau\le t\le (j+1)\tau \;\;
               {\rm and}\;\;j=0, 1, \ldots, m-1, \lb{3.9}\eea
as in p.114 of \cite{Lon1}. As in \cite{LoZ} and \cite{Lon3}, we denote the Maslov-type indices of
$\ga^m$ by $(i(\ga,m),\nu(\ga,m))$.

Then by Theorem 3.2 and some index computations of basic normal forms, the following iteration formula from \cite{LoZ} and \cite{Lon3} can be obtained.

\medskip

{\bf Theorem 3.3.} (cf. Theorem 9.3.1 of \cite{Lon3}) {\it For any path $\ga\in\mathcal{P}_\tau(2n)$,
let $M=\ga(\tau)$ and $C(M)=\sum_{0<\th<2\pi}S_M^-(e^{\sqrt{-1}\th})$. We extend $\ga$ to $[0,+\infty)$
by its iterates. Then for any $m\in\N$ we have
\bea &&i(\ga,m)
= m(i(\ga,1)+S^+_{M}(1)-C(M))\nn\\
&&\qquad+ 2\sum_{\th\in(0,2\pi)}E\left(\frac{m\th}{2\pi}\right)S^-_{M}(e^{\sqrt{-1}\th}) - (S_M^+(1)+C(M)) \lb{3.10}\eea
and
\be \hat{i}(\ga,1) = i(\ga,1) + S^+_{M}(1) - C(M) + \sum_{\th\in(0,2\pi)}\frac{\th}{\pi}S^-_{M}(e^{\sqrt{-1}\th}). \lb{3.11}\ee}

\medskip

{\bf Theorem 3.4.} (cf. Theorem 4.1 of \cite{LoZ} and Theorem 3.4 of \cite{DLLW}) {\it Fix an integer $q>0$. Let $\mu_i\ge 0$ and $\bb_i$ be integers for all $i=1,\cdots,q$. Let $\aa_{i,j}$
be positive numbers for $j=1,\cdots,\mu_i$ and $i=1,\cdots,q$. Let $\dl\in(0,\frac{1}{2})$ satisfying
$\dl\max\limits_{1\le i\le q}\mu_i<\frac{1}{2}$. Suppose $D_i \equiv \bb_i+\sum\limits_{j=1}^{\mu_i}\aa_{i,j}\neq 0$ for
$i=1,\cdots,q$. Then there exist infinitely many $(N, m_1,\cdots,m_q)\in\N^{q+1}$ such that
\bea
&& m_i\bb_i+\sum_{j=1}^{\mu_i}E(m_i\aa_{i,j}) =
      \varrho_i N+\Delta_i, \  \forall\ 1\le i\le q.  \lb{3.12}\\
&& \min\{\{m_i\aa_{i,j}\}, 1-\{m_i\aa_{i,j}\}\} < \dl,\ \forall\ j=1,\cdots,\mu_i, 1\le i\le q, \lb{3.13}\\
&& m_i\aa_{i,j}\in\N,\  {\rm if} \  \aa_{i,j}\in\Q,   \lb{3.14}\eea
where
\bea \varrho_i=\left\{\begin{array}{cc}1, &{\rm if}\ D_i>0, \cr
                                     -1, &{\rm if}\  D_i<0, \end{array}\right.\quad \Delta_i=\sum_{0<\{m_i\aa_{i,j}\}<\dl}1,\quad \forall\ 1\le i\le q.\lb{3.15}\eea}

{\bf Remark 3.5.} (i) When $D_i>0$ for all $1\le i\le q$, this is precisely the Theorem 4.1 of \cite{LoZ} (also cf.
Theorem 11.1.1 of \cite{Lon3}). When $D_i\neq0$ for all $1\le i\le q$, this has been proved in \cite{DLLW} (cf. Theorem 3.4 of \cite{DLLW}).

(ii) According to Theorem 3.4 and its proof, it is easy to see that for any two small enough integers $\dl_1$ and $\dl_2$ satisfying $\dl_k \in(0,\frac{1}{2})$ and $\dl_k\max\limits_{1\le i\le q}\mu_i<\frac{1}{2}$ for $k=1,2$, Theorem 3.4 holds with
$$\sum_{0<\{m_i\aa_{i,j}\}<\dl_1}1=\sum_{0<\{m_i\aa_{i,j}\}<\dl_2}1,\quad \forall\ 1\le i\le q.$$

\vspace{2mm}

In 2002, Y. Long and C. Zhu \cite{LoZ} has established the common index jump theorem for symplectic paths, which has become one
of the main tools to study the periodic orbit problem in Hamiltonian and symplectic dynamics. In \cite{DLW} of 2016, H. Duan,
Y. Long and W. Wang further improved this theorem to an enhanced version which gives more precise index properties of
$\ga_k^{2m_k}$ and $\ga_k^{2m_k\pm m}$ with $1\le m \le \bar{m}$ for any fixed $\bar{m}$. Under the help of Theorem 3.4, following
the proofs of Theorem 3.5 in \cite{DLW}, this result has been further generalized to the case of admitting the existence of
symplectic paths with negative mean indices.

\medskip

{\bf Theorem 3.6.} ({\bf Generalized common index jump theorem}, cf. Theorem 3.6 of \cite{DLLW})
{\it Let $\gamma_i\in\mathcal{P}_{\tau_i}(2n)$ for $i=1,\cdots,q$ be a finite collection of symplectic paths with nonzero mean
indices $\hat{i}(\ga_i,1)$. Let $M_i=\ga_i(\tau_i)$. We extend $\ga_i$ to $[0,+\infty)$ by (\ref{3.9}) inductively.

Then for any fixed $\bar{m}\in \N$, there exist infinitely many $(q+1)$-tuples
$(N, m_1,\cdots,m_q) \in \N^{q+1}$ such that the following hold for all $1\le i\le q$ and $1\le m\le \bar{m}$,
\bea
\nu(\ga_i,2m_i-m) &=& \nu(\ga_i,2m_i+m) = \nu(\ga_i, m),   \lb{3.16}\\
i(\ga_i,2m_i+m) &=& 2\varrho_i N+i(\ga_i,m),                         \lb{3.17}\\
i(\ga_i,2m_i-m) &=&  2\varrho_i N-i(\ga_i,m)-2(S^+_{M_i}(1)+Q_i(m)),  \lb{3.18}\\
i(\ga_i, 2m_i)&=& 2\varrho_i N -(S^+_{M_i}(1)+C(M_i)-2\Delta_i),     \lb{3.19}\eea
where \bea &&\varrho_i=\left\{\begin{array}{cc}1, &{\rm if}\ \hat{i}(\ga_i,1)>0, \cr
                                     -1, &{\rm if}\  \hat{i}(\ga_i,1)<0, \end{array}\right.\qquad
\Delta_i = \sum_{0<\{m_i\th/\pi\}<\delta}S^-_{M_i}(e^{\sqrt{-1}\th}),\nn\\
&&\ Q_i(m) = \sum_{e^{\sqrt{-1}\th}\in\sg(M_i),\atop \{\frac{m_i\th}{\pi}\}
                   = \{\frac{m\th}{2\pi}\}=0}S^-_{M_i}(e^{\sqrt{-1}\th}). \lb{3.20}\eea
Moreover we have
\bea \min\left\{\{\frac{m_i\theta}{\pi}\},1-\{\frac{m_i\theta}{\pi}\}\}\right\}<\dl,\lb{3.21}\eea
whenever $e^{\sqrt{-1}\theta}\in\sigma(M_i)$ and $\dl$ can be chosen as small as we want. More precisely, by (3.17) in \cite{DLLW} and (4.40), (4.41) in \cite{LoZ} , we have
\bea m_i=\left(\left[\frac{N}{M|\hat i(\gamma_i, 1)|}\right]+\chi_i\right)M,\quad\forall\  1\le i\le q,\lb{3.22}\eea
where $\chi_i=0$ or $1$ for $1\le i\le q$ and $\frac{M\theta}{\pi}\in\Z$
whenever $e^{\sqrt{-1}\theta}\in\sigma(M_i)$ and $\frac{\theta}{\pi}\in\Q$
for some $1\le i\le q$.  Furthermore, given $M_0$, from the proof of Theorem 4.1 of \cite{LoZ}, we may further require $N$ to be the mutiple of $M_0$, i.e., $M_0|N$.}

\vspace{2mm}

{\bf Remark 3.7.} In fact, let $\mu_i=\sum_{0<\th<2\pi}S_{M_i}^-(e^{\sqrt{-1}\th})$, $\alpha_{i,j}=\frac{\th_j}{\pi}$
where $e^{\sqrt{-1}\th_j}\in\sigma(M_i)$ for $1\le j\le\mu_i$ and $1\le i\le q$. Let $l=q+\sum_{i=1}^q \mu_i$ and
{\small\bea &&v=\left(\frac{1}{M|\hat{i}(\gamma_1,1)|},\cdots,\frac{1}{M|\hat{i}(\gamma_1,1)|},\right.\nn\\
&&\quad\quad\left.\frac{\aa_{1,1}}{|\hat{i}(\gamma_1,1)|},\cdots,
\frac{\aa_{1,\mu_1}}{|\hat{i}(\gamma_1,1)|},\cdots,\frac{\aa_{q,1}}{|\hat{i}(\gamma_q,1)|},\cdots,
\frac{\aa_{q,\mu_q}}{|\hat{i}(\gamma_q,1)|}\right)\in\R^l.\lb{3.23}\eea}
Then Theorem 3.6 is equivalent to find a vertex
\bes \chi=(\chi_1,\cdots,\chi_q,\chi_{1,1},\cdots,\chi_{1,\mu_1},\cdots,\chi_{q,1},\cdots,\chi_{q,\mu_q})\in\{0,1\}^l\ees
of the cube $[0,1]^l$ and infinitely many $N\in\N$ such that for any small enough $\epsilon\in(0,\frac{1}{2})$ there holds
\bea|\{Nv\}-\chi|<\epsilon.\lb{3.24}\eea

Next we use some multiple $\hat{N}$ of $N$ to replace $N$ in Theorem 3.6, then we get the corresponding iterations $\hat{m}_k$s such that the  equalities similar to (\ref{3.16})-(\ref{3.19}) hold. At the same time, we also get some relation equalities about these integers, which will play a crucial role in the proof of Theorem 1.1 in Section 4.

\vspace{2mm}

{\bf Theorem 3.8.} (A variant of Theorem 3.6)  {\it Let $\gamma_i\in\mathcal{P}_{\tau_i}(2n)$ for $i=1,\cdots,q$ be a finite collection of symplectic paths with nonzero mean
indices $\hat{i}(\ga_i,1)$. Let $M_i=\ga_i(\tau_i)$. Fix a positive integer $\hat{p}$, use $\frac{\delta}{\hat{p}}$ and $\frac{\epsilon}{\hat{p}}$ to substitute $\delta$ in (\ref{3.21}) and $\epsilon$ in (\ref{3.24}) respectively, and then let $N$ be an integer satisfying (\ref{3.16})-(\ref{3.22}).

Then for $\hat{N}=\hat{p}N$, there exists the corresponding $q$-tuples $(\hat{m}_1,\cdots,\hat{m}_q) \in \N^q$ such that the following hold for all $1\le i\le q$ and $1\le m\le \bar{m}$
\bea
\nu(\ga_i,2\hat{m}_i-m) &=& \nu(\ga_i,2\hat{m}_i+m) = \nu(\ga_i, m),   \lb{3.25}\\
i(\ga_i,2\hat{m}_i+m) &=& 2\varrho_i \hat{p} N+i(\ga_i,m),                         \lb{3.26}\\
i(\ga_i,2\hat{m}_i-m) &=&  2\varrho_i \hat{p} N-i(\ga_i,m)-2(S^+_{M_i}(1)+Q_i(m)),  \lb{3.27}\\
i(\ga_i, 2\hat{m}_i)&=& 2\varrho_i \hat{p} N -(S^+_{M_i}(1)+C(M_i)-2\hat{\Delta}_i),     \lb{3.28}\eea
where $\varrho_i$ and $Q_i(m)$ are the same as those in (\ref{3.20}) and for any $1\le i\le q$ \bea \hat{m}_i=\left(\left[\frac{\hat{p}N}{M|\hat i(\gamma_i, 1)|}\right]+\hat{\chi}_i\right)M,\quad
\hat{\Delta}_i = \sum_{0<\{\hat{m}_i\th/\pi\}<\delta}S^-_{M_i}(e^{\sqrt{-1}\th}).\lb{3.29}\eea
Furthermore, comparing with some corresponding integers in Theorem 3.6, there holds
\bea \hat{\chi}_i=\chi_i,\quad \hat{m}_i=\hat{p}m_i,\quad \hat{\Delta}_i=\Delta_i,\quad \forall\ 1\le i\le q.\lb{3.30}\eea}

{\bf Proof.} Set $v_i=\frac{1}{|M\hat{i}(\gamma_i,1)|},\forall 1\le i\le q$. Firstly notice that \bea\{\hat{p}Nv_i\}=\{\hat{p}([Nv_i]+\{Nv_i\})\}=\{\hat{p}[Nv_i]+\hat{p}\{Nv_i\}\}=\{\hat{p}\{Nv_i\}\}.\eea

When $\chi_i=0$, it follows from (\ref{3.24}) with $\epsilon$ replaced by $\frac{\epsilon}{\hat{p}}$ (i.e. $\{Nv_i\}<\frac{\epsilon}{\hat{p}}$) that  $|\{\hat{p}Nv_i\}-\chi_i|=\{\hat{p}\{Nv_i\}\}=\hat{p}\{Nv_i\}<\epsilon$. When $\chi_i=1$, it follows from (\ref{3.24}) that $0<1-\{Nv_i\}<\frac{\epsilon}{\hat{p}}$. Therefore there holds $0<\hat{p}-\hat{p}\{Nv_i\}<\epsilon$, which yields $0<1-\{\hat{p}\{Nv_i\}\}=1-\{\hat{p}Nv_i\}<\epsilon$. In a word, there always holds
\bea |\{\hat{p}Nv_i\}-\chi_i|=|\{\hat{p}\{Nv_i\}\}-\chi_i|<\epsilon,\qquad 1\le i\le q.\lb{3.32}\eea
Roughly speaking, when $\{Nv_i\}$ is close enough to $\chi_i$, $\{\hat{p}Nv_i\}$ is also close enough to the same $\chi_i$. Therefore $\hat{\chi}_i=\chi_i, \forall\ 1\le i\le q$.

Now we prove the relation of $\hat{m}_i$ and $m_i$ for $1\le i\le q$. By the similar arguments as those in proof of Theorem 3.6 (i.e. the proof of similar (\ref{3.22})), we have
\bea \hat{m}_i&=&\left(\left[\frac{\hat{p}N}{M|\hat i(\gamma_i, 1)|}\right]+\hat{\chi}_i\right)M=\left(\left[\hat{p}Nv_i\right]+\chi_i\right)M\nn\\
  &=&\left(\left[\hat{p}([Nv_i]+\{Nv_i\})\right]+\chi_i\right)M\nn\\
  &=&\left(\hat{p}[Nv_i]+\left[\hat{p}\{Nv_i\}\right]+\chi_i\right)M.\lb{3.33}\eea

When $\chi_i=0$, then $\{Nv_i\}<\frac{\epsilon}{\hat{p}}$, and thus $\hat{p}\{Nv_i\}<\epsilon$. By (\ref{3.33}) and (\ref{3.22}) we obtain $\hat{m}_i=\hat{p}[Nv_i]M=\hat{p}m_i$. When $\chi_i=1$, it follows from (\ref{3.24}) that $0<1-\{Nv_i\}<\frac{\epsilon}{\hat{p}}$. Therefore there holds $\hat{p}-1<\hat{p}-\epsilon<\hat{p}\{Nv_i\}<\hat{p}$, which yields $[\hat{p}\{Nv_i\}]=\hat{p}-1$. Again by (\ref{3.33}) and (\ref{3.22}) we obtain $\hat{m}_i=(\hat{p}[Nv_i]+(\hat{p}-1)+1)M=\hat{p}([Nv_i]+1)M=\hat{p}([Nv_i]+\chi_i)M=\hat{p}m_i$. In a word, there always holds $\hat{m}_i=\hat{p}m_i$, $\forall 1\le i\le q$.

Finally the proof of identities (\ref{3.25})-(3.28) are completely the same as those of (\ref{3.16})-(\ref{3.19}) in Theorem 3.6, so here we omit these details. \hfill\hb

\setcounter{figure}{0}
\setcounter{equation}{0}
\section{Proof of Theorem 1.1}

In this section, let $F$ be any bumpy Finsler metric on $S^3$, i.e., all closed geodesics (including iterations) on $(S^3,F)$ are non-degenerate. To prove theorem 1.1, firstly we assume that the Morse index $i(c)\ge 1$ for any prime closed geodesic $c$ on $(S^3,F)$.

On one hand, Theorem 1.4 in \cite{DLW} showed that if every prime closed geodesic $c$ on every bumpy $S^3$ satisfies $i(c)\ge 2$, then there exist at least four prime closed geodesics. On the other hand, Remark 3.7 in \cite{DuL2} showed that if there exist exactly two prime closed geodesics $c_1$ and $c_2$ on bumpy Finsler $S^3$, then only two possible cases happen: (i) $i(c_1)=0$ and $i(c_2)=1$; (ii) $i(c_k)\ge 2$ for $k=1,2$.

Thus, it follows from Remark 3.7 in \cite{DuL2} and the assumption of nonzero Morse index that there exist at least 3 prime closed geodesics on $S^3$. Therefore, in order to prove theorem 1.1, by the above assumption and Theorem 1.4 of \cite{DLW}, we only need to consider the case when there exist exactly 3 prime closed geodesics $c_k$ with $k=1,2,3$ on $(S^3,F)$ and at least one prime closed geodesic among them has Morse index $1$. Then we try to get a contradiction.

Furthermore, by Propositions 2.1, 2.3 and 2.4, among three prime closed geodesics there exists at least one prime closed geodesic with even Morse index such that these closed geodesics can generate non-trivial Morse-type numbers $M_{2j}\ge b_{2j}\ge 1$ for any $j\ge 1$. In fact, under the assumption of existence of exactly three prime closed geodesics, it can be further showed that there exist exactly two prime closed geodesics with even Morse indices (cf. Lemma 4.2 and Assumption below).

It is well-known that the Morse index sequence $\{i(c^m)\}$ either tends to $+\infty$ asumptotically linearly or $i(c^m)=0$ for all $m\ge1$ (cf. Proposition 1.3 of \cite{Bot}). Therefore $\hat{i}(c)=0$ if and only if $i(c^m)=0$ for all $m\ge1$. Since $i(c_k)\ge 1$ for $k=1,2,3$ by the above assumption, so the mean index satisfies $\hat{i}(c_k)>0$. So $\varrho_k=1$ in Theorem 3.6 by the definition (\ref{3.20}), and $i(c_k^m)\rightarrow +\infty$ as $m\rightarrow +\infty$. So the positive integer $\hat{m}$ defined by
\bea \overline{m}=\max\limits_{1\leq k\leq q}\min\{m_0\in\N\,|\,i(c_{k}^{m+m_0})\geq i(c_{k})+4,\  \forall\, m\ge 1\}\lb{4.1}\eea
is well defined and finite.

{\bf Lemma 4.1.} {\it With the above $\overline{m}$, by Theorem 3.6 there exist $N\in\N$ and $q$-tuples $(m_1,\cdots,m_q) \in \N^{q}$ such that the following inequalities hold}
\bea i(c_{k}^{2m_{k}-m})&\leq& 2N-i(c_{k}),\quad\forall\ 1\leq m< 2m_{k}, \lb{4.2}\\
 i(c_{k}^{2m_{k}+m})&\geq& 2N+i(c_{k}),\quad \forall\ m\geq 1.\lb{4.3}\eea

{\bf Proof.} Note that there always holds $S^+_{M_i}(1)=Q_i(m)=0$ by the non-degeneracy of $F$ and Lemma 3.1, and $i(c^m)\ge i(c)\ge 1$ for any $m\ge 1$ by the Bott formulae (cf. Theorem 9.2.1 of \cite{Lon3}). So it follows from (\ref{3.17}) and (\ref{3.18}) of Theorem 3.6 that the above inequalities hold for $m\leq \overline{m}$.

Now assume $m\geq \overline{m}+1$. Similar to the proof of Theorem 3.6 (cf. equalities (3.33) and (3.36) in \cite{DLW} for more details), by Theorem 3.3 and (3.12), we have
\bea i(c_{k}^{2m_{k}-m})&=&-i(c_{k}^{m})+2N+2\Delta_{k}+2\sum\limits_{\theta}\left(E\left(\frac{2m_{k}-m}{2\pi}\theta\right)+\right.\nn\\ &&\quad\left. E\left(\frac{m}{2\pi}\theta\right)-E\left(\frac{m_{k}}{\pi}\theta\right)\right)S_{M_{k}}^{-}(e^{\sqrt{-1}\theta})-2C(M_{k}),\lb{4.4}\\
i(c_{k}^{2m_{k}+m})&=&i(c_{k}^{m})+2N+2\Delta_{k}+2\sum\limits_{\theta}\left(E\left(\frac{2m_{k}+m}{2\pi}\theta\right)-\right.\nn\\ &&\qquad \left. E\left(\frac{m}{2\pi}\theta\right)-E\left(\frac{m_{k}}{\pi}\theta\right)\right)S_{M_{k}}^{-}(e^{\sqrt{-1}\theta}).\lb{4.5}\eea

Note that $E(a)+E(b)-E(a+b)\leq 1$. Then by the definition of $\overline{m}$, there holds
\bea i(c_{k}^{2m_{k}-m})&\le& -i(c_{k}^{m})+2N+2\Delta_{k}\nn\\
&\le& -i(c_{k})-4+2N+2\Delta_{k}\leq 2N-i(c_{k}),\\
i(c_{k}^{2m_{k}+m})&\ge& i(c_{k}^{m})+2N+2\Delta_{k}-2C(M_{k})\nn\\
&\ge& i(c_{k})+4+2N+2\Delta_{k}-2C(M_{k})\geq 2N+i(c_{k}),\eea
where we use the fact $0\le\Delta_{k}\le C(M_k)\le 2$ in the case of $S^3$. \hfill\hb

\vspace{2mm}

{\bf Lemma 4.2.} {\it It is impossible that two of $c_{1},c_{2},c_{3}$ have odd indices and one of them has an even index.}

\vspace{2mm}

{\bf Proof.} Without loss of generality, suppose that $c_{1},c_{2}$ have odd indices and $c_{3}$ has an even index, then by Lemma 4.1, there exist positive integers $N$ and $m_{3}$ such that
\bea i(c_{3}^{2m_{k}-m})&\le& 2N-i(c_{3})\le 2N-1,\quad \forall\  1\leq m< 2m_{3}, \lb{4.8}\\
 i(c_{3}^{2m_{k}+m})&\ge& 2N+i(c_{3})\ge 2N+1,\quad \forall\  m\geq 1.\lb{4.9}\eea
Then by Proposition 2.1, $c_{1}$ and $c_{2}$ have no contribution to $M_{2N}$ while $c_{3}$ contributes at most 1 to $M_{2N}$. Hence there holds $M_{2N}\leq 1<b_{2N}$, which contradicts to Proposition 2.4.\hfill\hb

\vspace{2mm}

According to the above arguments and Lemma 4.2, we make the following assumption and will get a contradiction to complete the proof of Theorem 1.1.

\vspace{2mm}

{\bf Assumption.} {\it There exist exactly 3 prime closed geodesics $c_{1},c_{2},c_{3}$ on bumpy Finsler 3-sphere $(S^3,F)$ with $i(c_{1})=1, i(c_{2})\in2\N, i(c_{3})\in2\N$.}

\vspace{2mm}

Now we use Lemma 4.1 to get the following two crucial inequalities.

\vspace{2mm}

{\bf Lemma 4.3.} {\it Under the assumption, with the above $\overline{m}$, Theorem 3.6 yields $(q+1)$-tuples $(N, m_1,\cdots,m_q) \in \N^{q+1}$ such that the following inequalities hold
\bea \sum_{k=1}^{3}2m_{k}\gamma_{c_k}-(2N)_{+}^{e}+(2N)_{+}^{o}&\ge& 2N-1, \lb{4.10}\\
\sum_{k=1}^{3}2m_{k}\gamma_{c_k}-(2N-1)_{+}^{e}+(2N-1)_{+}^{o}&\le& 2N-3, \lb{4.11}\eea
where the counting notations are defined by $n_{+}^{e}=\#\{k\ |\ i(c_{k}^{2m_{k}})>n,i(c_{k}^{2m_{k}})\equiv i(c_{k})\equiv 0\left(\rm mod 2\it\right)\}$ and $n_{+}^{o}=\#\{k\ |\ i(c_{k}^{2m_{k}})>n,i(c_{k}^{2m_{k}})\equiv i(c_{k})\equiv 1\left(\rm mod 2\it\right)\}$.}

\vspace{2mm}

{\bf Proof.} It follows from Proposition 2.1, Definition 2.2 and Theorem 3.3 that
\bea \sum_{m=1}^{2m_k} (-1)^{i(c_k^m)} \dim \ol{C}_{i(c_k^{m})}(E,c_k^m)
&=& \sum_{i=0}^{m_k-1} \sum_{m=2i+1}^{2i+2} (-1)^{i(c_k^m)} \dim \ol{C}_{i(c_k^{m})}(E,c_k^m) \nn\\
&=& \sum_{i=0}^{m_k-1} \sum_{m=1}^{2} (-1)^{i(c_k^m)} \dim \ol{C}_{i(c_k^{m})}(E,c_k^m) \nn\\
&=& m_k \sum_{m=1}^{2} (-1)^{i(c_k^m)} \dim \ol{C}_{i(c_k^{m})}(E,c_k^m) \nn\\
&=& 2m_k\ga_{c_k},\qquad \forall\ 1\le k\le 3, \lb{4.12}\eea
where the second equality follows from Proposition 2.1 and the fact $i(c_k^{m+2})-i(c_k^m)\in 2\Z$ for all
$m\ge 1$ from Theorem 3.3, and the last equality follows from Proposition 2.1 and Definition 2.2.

By Lemma 4.1 and Proposition 2.1, the sum of the left side of (\ref{4.12}) with respect to $k$ is
\bea &&\sum_{k=1}^{3}\sum_{m=1}^{2m_{k}}(-1)^{i(c_{k}^{m})}\dim\overline{C}_{i(c_k^{m})}(E,c_{k}^{m})\nn\\
&&\qquad=\sum_{i=1}^{2N}(-1)^{i}M_{i}+(2N)_{+}^{e}-(2N)_{+}^{o},\lb{4.13}\\
&&\sum_{k=1}^{3}\sum_{m=1}^{2m_{k}}(-1)^{i(c_{k}^{m})}\dim\overline{C}_{i(c_k^{m})}(E,c_{k}^{m})\nn\\
&&\qquad=\sum_{i=1}^{2N-1}(-1)^{i}M_{i}+(2N-1)_{+}^{e}-(2N-1)_{+}^{o}.\lb{4.14}\eea

On the other hand, Propositions 2.3 and 2.4 give the following inequalities
\bea \sum_{i=1}^{2N}(-1)^{i}M_{i}\ge \sum_{i=1}^{2N}(-1)^{i}b_{i}=2N-1,\lb{4.15}\\
\sum_{i=1}^{2N-1}(-1)^{i}M_{i}\le \sum_{i=1}^{2N-1}(-1)^{i}b_{i}=2N-3.\lb{4.16}\eea
Now the inequalities (\ref{4.10}) and (\ref{4.11}) follow from (\ref{4.12})-(\ref{4.16}). \hfill\hb

\vspace{2mm}

Again by Lemma 4.1 and Proposition 2.1, $c_{2}^{2m_{2}}$ and $c_{3}^{2m_{3}}$ contribute at most 1 to $M_{2N}$, respectively. Note that Propositions 2.3 and 2.4 yields $M_{2N}\geq b_{2N}=2$. So $c_{2}^{2m_{2}}$ and $c_{3}^{2m_{3}}$ must contribute exactly 1 to $M_{2N}$, respectively, that is, there must hold
\bea i(c_{2}^{2m_{2}})=2N,\qquad i(c_{3}^{2m_{3}})=2N. \lb{4.17}\eea

Hence, by the definitions of $n_{+}^{e}$ and $n_{+}^{o}$  we have $(2N)_{+}^{e}=0$, $(2N-1)_{+}^{e}=2$, $(2N)_{+}^{o}\le 1$ and $0\le (2N-1)_{+}^{o}$. And then it follows from (\ref{4.10}) and (\ref{4.11}) of Lemma 4.3 that
\bea 2N-2&\leq&2N-1-(2N)_{+}^{o}\nonumber\\
&\leq&\sum_{k=1}^{3}2m_{k}\gamma_{c_k}\nonumber\\
&\leq& 2N-1-(2N-1)_{+}^{o}\nonumber\\
&\leq&2N-1.\lb{4.18}\eea

Fix a small enough $\epsilon$ in (3.24) and suppose $N$ satisfies (3.24) for $\frac{\epsilon}{4}$ (i.e., taking $\hat{p}=4$ there), then by Theorem 3.8, the integers $\hat{N}=4N$ and $\hat{m}_{k}=\left(\left[\frac{\hat{N}}{M\hat{i}(c_k)}\right]+\chi_{k}\right)M=4m_k, k=1,\cdots,q$ yields the $(q+1)$-tuples $(\hat{N},\hat{m}_1,\cdots,\hat{m}_q) \in \N^{q+1}$ such that Theorem 3.4 and (\ref{3.25})-(\ref{3.28}) in Theorem 3.8 hold. Then by the similar arguments above, it can be shown that Lemma 4.1, Lemma 4.3 and (\ref{4.18}) also hold with $N$ being replaced by $4N$. This is to say, we have
\bea 8N-2\leq \sum_{k=1}^{3}2\hat{m}_{k}\gamma_{c_k}=4\sum_{k=1}^{3}2 m_{k}\gamma_{c_k}\leq 8N-1,\lb{4.19}\eea where we use the fact $\hat{m}_{k}=4m_{k}$ in (\ref{3.30}) of Theorem 3.8. Then the integer in the middle is the multiple of $4$ and lies between $8N-2$ and $8N-1$. Clearly this is a contradiction.

Therefore the Assumption is not true and the proof of Theorem 1.1 is finished. \hfill\hb

\vspace{1mm}

{\bf Acknowledgements.} The first author would like to thank Professor Yiming Long and Chongchun Zeng for their valuable communications and constant encouragement.






\bibliographystyle{abbrv}

\end{document}